\documentclass[12pt,reqno]{amsart}
\usepackage{fullpage,mathpple}

\newcommand{\ep}{\varepsilon}
\newcommand{\li}{\mathop{\rm li}}
\renewcommand{\P}[1]{\mathcal P_{#1}}
\newcommand{\N}{{\mathbb N}}

\newcommand{\Z}{{\mathbb Z}}
\renewcommand{\mod}[1]{{\ifmmode\text{\rm\
(mod~$#1$)}\else\discretionary{}{}{\hbox{ }}\rm(mod~$#1$)\fi}}
\newsymbol\dnd 232D
\newcommand{\lcm}{\mathop{\rm{lcm}}}
\newcommand{\ord}{\operatorname{ord}}
\newcommand{\oo}{\operatorname{ord^*}}

\newtheorem{theorem}{Theorem}
\newtheorem{conj}[theorem]{Conjecture}
\newtheorem{lemma}[theorem]{Lemma}
\newtheorem{prop}[theorem]{Proposition}

\begin{document}

\title
[The iterated Carmichael $\lambda$-function and cycles of the power  
generator]
{The iterated
Carmichael~$\lambda$-function
and the number of cycles of the power generator}
\author{Greg Martin}
\address{Department of Mathematics\\University of British
Columbia\\Room 121, 1984 Mathematics Road\\Vancouver, BC V6T 1Z2\\
Canada}
\email{gerg@math.ubc.ca}

\author{Carl Pomerance}
\address{Mathematics Department\\
Dartmouth College\\
Hanover, NH 03755-3551\\
U.S.A.}
\email{carlp@math.dartmouth.edu}

\thanks{G.M.\ is supported in part by the National Sciences and Engineering Research Council of Canada. C.P.\ is supported in part by the National Science Foundation.}

\maketitle

\section{introduction}

A common pseudorandom number generator is the power generator:
$x\mapsto x^\ell\mod n$.  Here, $\ell,n$ are fixed integers at
least 2, and one constructs a pseudorandom sequence by starting
at some residue mod~$n$ and iterating this $\ell$th power map.
(Because it is the easiest to compute, one often takes $\ell=2$;
this case is known as the BBS generator, for Blum, Blum, and Shub.)
To be a good generator, the period should be large.  Of course,
the period depends somewhat on the number chosen for the initial value.  
  However,
a universal upper bound for this period is $\lambda(\lambda(n))$ where
$\lambda$ is Carmichael's function.  Here, $\lambda(m)$ is defined as  
the
order of the largest cyclic subgroup of the multiplicative group
$(\Z/m\Z)^\times$.  It may be computed via the identity
$\lambda(\lcm\{a,b\})=\lcm\{\lambda(a),\lambda(b)\}$ and its values
at prime powers: with $\phi$ being Euler's function, 
$\lambda(p^a)=\phi(p^a)=(p-1)p^{a-1}$ for every odd prime power $p^a$  
and for 2 and 4,
and $\lambda(2^a)=\phi(2^a)/2=2^{a-2}$ for $a\ge 3$.

Statistical properties of $\lambda(n)$ were studied by Erd\H os,  
Schmutz, and the second author in \cite{EPS}, and
in particular, they showed that
$\lambda(n)=n/\exp((1+o(1))\log\log n\log\log\log n)$ as $n\to\infty$
through a certain set of integers of asymptotic density 1.  This does
not quite pinpoint the normal order of $\lambda(n)$ (even the sharper version of this
theorem from \cite{EPS} falls short in this regard), but it is certainly
a step in this direction, and does give the normal order of the function
$\log(n/\lambda(n))$.

In this paper we prove a result of similar quality for the function
$\lambda(\lambda(n))$, which we have seen arises in connection with
the period of the power generator.  We obtain the same expression as with $\lambda(n)$,
except that the $\log\log n$ is squared.
That is, $\lambda(\lambda(n))=n/\exp((1+o(1))(\log\log n)^2\log\log\log  
n)$
almost always.

We are able to use this result to say something nontrivial about the
number of cycles for the power generator.  This problem has been  
considered
in several papers, including \cite{BHM}, \cite{BG}, and \cite{R}.
We show that for almost all integers $n$, the number of cycles for the
$\ell$th power map modulo $n$ is at least 
$\exp((1+o(1))(\log\log n)^2\log\log\log n)$,
and we conjecture that this lower bound is actually the truth.  Under the assumption of
the Generalized Riemann Hypothesis (GRH), and using a new result of  
Kurlberg and the second author \cite{KP}, we prove our conjecture.
(By the GRH, we mean the Riemann Hypothesis for Kummerian fields
as used by Hooley in his celebrated conditional proof of the Artin
conjecture.)

For an arithmetic function $f(n)$ whose values are in the natural  
numbers,
let $f_k(n)$ denote the $k$th iterate of $f$ evaluated at $n$.
One might ask about the normal behavior of $\lambda_k(n)$ for $k\ge 3$.
Here we make a conjecture for each fixed $k$.
We also briefly consider the function $L(n)$ defined as the least $k$  
such
that $\lambda_k(n)=1$.  A similar undertaking was made by Erd\H os,  
Granville, Spiro, and the second author in \cite{EGPS}
for the function $F(n)$ defined as the least $k$ with $\phi_k(n)=1$.
Though $\lambda$ is very similar to $\phi$, the behavior
of $L(n)$ and $F(n)$ seem markedly different.  We know that $F(n)$
is always of order of magnitude $\log n$, and it is shown in  
\cite{EGPS},
assuming the Elliott--Halberstam conjecture on the average distribution  
of
primes in arithmetic progressions with large moduli, that in fact  
$F(n)\sim\alpha\log n$ on a set of asymptotic density 1 for a  
particular positive constant $\alpha$.
We know far less about $L(n)$, not even its typical order of magnitude.
We raise the possibility that it is normally of order $\log\log n$
and show that it is bounded by this order infinitely often.

A more formal statement of our results follows.

\begin{theorem}
\label{lambda.lambda.normal.order.theorem}
The normal order of $\log\big( n/\lambda(\lambda(n)) \big)$ is
$(\log\log n)^2 \log\log\log n$.
That is,
\[
\lambda(\lambda(n))=n\exp\left(-(1+o(1))(\log\log n)^2\log\log\log  
n\right)
\]
as $n\to\infty$ through a set of integers of asymptotic density $1$.
\end{theorem}
We actually prove the slightly stronger result:
given any function
$\psi(n)$ going to infinity arbitrarily slowly, we have
\[
\lambda(\lambda(n)) = n\exp\!\big( {-(\log\log n)^2} ( \log\log\log n +
O(\psi(n)) ) \big)
\]
for almost all $n$.

Given integers $\ell,n\ge2$, let $C(\ell,n)$ denote the number of cycles
when iterating the modular power map $x\mapsto x^\ell\mod n$.

\begin{theorem}
\label{number.of.cycles.theorem}
Given any fixed integer $\ell\ge2$, there is a set of integers of  
asymptotic
density $1$ such that as $n\to\infty$ through this set,
\begin{equation}
\label{numberofcycles}
C(\ell,n) \ge \exp\!\big((1+o(1)) (\log\log n)^2 \log\log\log n  \big).
\end{equation}
Further, if $\ep(n)$ tends to $0$ arbitrarily slowly, we have
$C(\ell,n)\le n^{1/2-\ep(n)}$ for almost all $n$.  Moreover,
for a positive proportion of integers $n$ we have $C(\ell,n)\le  
n^{.409}$.
Finally, if the Generalized Riemann Hypothesis (GRH) is true,
we have equality in \eqref{numberofcycles} on a set of integers $n$ of  
asymptotic density $1$.
\end{theorem}

\begin{conj}
\label{lambdak.normal.order.conj}
The normal order of $\log(n/\lambda_k(n))$ is $(1/(k-1)!)(\log\log
n)^k\log\log\log n$. That is, for each fixed integer $k\ge1$,
\[
\lambda_k(n) =
n\exp\!\left({-\left(\frac1{(k-1)!}+o(1)\right)(\log\log  
n)^k}(\log\log\log n)\right)
\]
for almost all $n$.
\end{conj}

Define $L(n)$ to be the number of iterations of $\lambda$ required to
take $n$ to 1, that is, $L(n)$ equals the smallest nonnegative integer
$k$ such that $\lambda_k(n)=1$.

\begin{theorem}
\label{log.log.iterates.theorem}
There are infinitely many integers $n$ such that $L(n) < (1/\log  
2+o(1))\log\log n$.
\end{theorem}

\section{Notation, strategy, and preliminaries}

The proof of Theorem \ref{lambda.lambda.normal.order.theorem}, our  
principal result,
proceeds by comparing the prime divisors of
$\lambda(\lambda(n))$ with those of $\phi(\phi(n))$. The primes
dividing $\phi(m)$ and $\lambda(m)$ are always the same. However, this
is not always true for $\phi(\phi(m))$ and $\lambda(\lambda(m))$. The
prime 2 clearly causes problems; for example, we have $\phi(\phi(8))=2$
but $\lambda(\lambda(8))=1$. However this problem also arises from the
interaction between different primes, for example, $\phi(\phi(91))=24$
but $\lambda(\lambda(91))=2$.

We shall use the following notation throughout the paper.
The letters $p,q,r$ will always denote primes.  Let $v_q(n)$ denote the
exponent on $q$
in the prime factorization of $n$, so that
\[
n~=~\prod_q q^{v_q(n)}
\]
for every positive integer $n$. We let $\P n = \{p\colon p\equiv1\mod
n\}$. We let $x>e^{e^e}$ be a real number and $y=y(x)=\log\log x$.  By
$\psi(x)$ we denote a function tending to infinity but more slowly than
$\log\log\log x=\log y$. In Sections 2--5, the phrase ``for almost all  
$n$'' always means ``for all but $O(x/\psi(x))$ integers $n\le x$''.

First we argue that the ``large'' prime divisors typically do not
contribute significantly:

\begin{prop}
For almost all $n\le x$, the prime divisors of $\phi(\phi(n))$ and
$\lambda(\lambda(n))$ that exceed $y^2$ are identical.
\label{same.large.prime.divisors.prop}
\end{prop}

\begin{prop}
For almost all $n\le x$,
\begin{equation}
\sum_{\substack{q>y^2 \\ v_q(\phi(\phi(n)))\ge2}} v_q(\phi(\phi(n)))
\log q~\ll~y^2\psi(x).
\label{large.primes.in.phi.phi.eq}
\end{equation}
\label{large.primes.in.phi.phi.prop}
\end{prop}

Next we argue that the contribution of ``small'' primes to
$\lambda(\lambda(n))$ is typically small:

\begin{prop}
For almost all $n\le x$, we have
\[
\sum_{q\le y^2} v_q(\lambda(\lambda(n))) \log q~\ll~y^2\psi(x).
\]
\label{small.primes.in.lambda.lambda.prop}
\end{prop}

Finally, we develop an understanding of the typical contribution of 
small primes
to $\phi(\phi(n))$ by comparing it to the additive function
$h(n)$ defined by
\begin{equation}
h(n)~=~\sum_{p\mid n} \sum_{r\mid p-1} \sum_{q\le y^2} v_q(r-1) \log q.
\label{h.definition}
\end{equation}

\begin{prop}
For almost all $n\le x$,
\[
\sum_{q\le y^2} v_q(\phi(\phi(n))) \log q ~=~h(n)+O(y\log  
y\cdot\psi(x)).
\]
\label{h.good.approximation.prop}
\end{prop}

\begin{prop}
For almost all $n\le x$, we have $h(n) = y^2\log y + O(y^2)$.
\label{h.normal.order.prop}
\end{prop}

\begin{proof}[Proof of Theorem \ref{lambda.lambda.normal.order.theorem}]
Let $x$ be a sufficiently large real number. For any positive integer  
$n\le x$
we may write
\begin{equation*}
\log \frac n{\lambda(\lambda(n))}~=~\log \frac n{\phi(n)} + \log
\frac{\phi(n)}{\phi(\phi(n))} + \log
\frac{\phi(\phi(n))}{\lambda(\lambda(n))}.
\end{equation*}
Recall that $n/\phi(n) \ll \log\log n$, and so the first two
terms are both $O(\log\log\log x)$.  Thus, it suffices to show that
\begin{equation}
\log \frac{\phi(\phi(n))}{\lambda(\lambda(n))}~=~(\log\log x)^2
(\log\log\log x+O(\psi(x)))~=~y^2\log y + O(y^2\psi(x))
\label{goal.equation}
\end{equation}
for almost all $n\le x$. We write
\begin{align}
\log\frac{\phi(\phi(n))}{\lambda(\lambda(n))}~ &=~\sum_{q} \big(
v_q(\phi(\phi(n))) - v_q(\lambda(\lambda(n))) \big) \log q \notag \\
&=~\sum_{q\le y^2} v_q(\phi(\phi(n))) \log q - \sum_{q\le y^2}
v_q(\lambda(\lambda(n))) \log q
\label{split.into.large.and.small.primes} \\
&\qquad+ \sum_{q>y^2} \big( v_q(\phi(\phi(n))) -
v_q(\lambda(\lambda(n))) \big) \log q. \notag
\end{align}
Since $\lambda(\lambda(n))$ always divides $\phi(\phi(n))$, the
coefficients of $\log q$ in this last sum are all nonnegative. On the
other hand, Proposition \ref{same.large.prime.divisors.prop} tells us
that for almost all $n\le x$, whenever $v_q(\phi(\phi(n))) > 0$ we have
$v_q(\lambda(\lambda(n))) > 0$ as well. Therefore the primes $q$ for
which $v_q(\phi(\phi(n))) \le 1$ do not contribute to this last sum at
all, that is,
\begin{align*}
0~&\le~\sum_{q>y^2} \big( v_q(\phi(\phi(n))) - v_q(\lambda(\lambda(n)))
\big) \log q \\
&=~\sum_{\substack{q>y^2 \\ v_q(\phi(\phi(n)))\ge2}} \big(
v_q(\phi(\phi(n))) - v_q(\lambda(\lambda(n))) \big) \log q \\
&\le~\sum_{\substack{q>y^2 \\ v_q(\phi(\phi(n)))\ge2}}
v_q(\phi(\phi(n))) \log q \ll y^2\psi(x)
\end{align*}
for almost all $n\le x$ by Propositions \ref{same.large.prime.divisors.prop} and
\ref{large.primes.in.phi.phi.prop}. Moreover, Proposition
\ref{small.primes.in.lambda.lambda.prop} tells us that the second sum
on the right-hand side of equation
\eqref{split.into.large.and.small.primes} is $O(y^2\psi(x))$ for almost
all $n\le x$. Therefore equation
\eqref{split.into.large.and.small.primes} becomes
\[
\log\frac{\phi(\phi(n))}{\lambda(\lambda(n))}~=~\sum_{q\le y^2}
v_q(\phi(\phi(n))) \log q + O(y^2\psi(x))
\]
for almost all $n\le x$. By Proposition
\ref{h.good.approximation.prop}, the sum on the right-hand side can be
replaced by $h(n)$ for almost all $n\le x$, the error $O(y\log y\cdot
\psi(x))$ in that proposition being absorbed into the existing error
$O(y^2\psi(x))$. Finally, Proposition \ref{h.normal.order.prop} tells
us that $h(n) = y^2\log y + O(y^2)$ for almost all $n\le x$. We
conclude that equation \eqref{goal.equation} is satisfied for almost
all $n\le x$, which establishes the theorem.
\end{proof}

Given integers $a$ and $n$, recall that $\pi(t;n,a)$ denotes the number
of primes up to $t$ that are congruent to $a\mod n$. The
Brun--Titchmarsh inequality (see \cite[Theorem 3.7]{HR}) states that
\begin{equation}
\pi(t;n,a)~\ll~\frac t{\phi(n)\log(t/n)}
\label{real.BT}
\end{equation}
for all $t>n$. We use repeatedly a weak form of this inequality, valid
for all $t>e^e$,
\begin{equation}
\sum_{\substack{p\le t\\ p\in\P n}}\frac1p~\ll~\frac{\log\log  
t}{\phi(n)},
\label{primesum}
\end{equation}
which follows from the estimate \eqref{real.BT} with $a=1$ by partial
summation. When $n/\phi(n)$ is bounded, this estimate simplifies to
\begin{equation}
\sum_{\substack{p\le t\\ p\in\P n}} \frac1p ~\ll ~\frac {\log\log t}n.
\label{BT.prime.power}
\end{equation}
For example, we shall employ this last estimate when $n$ is a prime or
a prime power and when $n$ is the product of two primes or prime
powers; in these cases we have $n/\phi(n)\le3$. We also quote the fact
(see Norton \cite{N} or the paper \cite{P} of the second author) that
\begin{equation}
\sum_{\substack{p\in\P n \\ p\le t}} \frac1{p}~ =~ \frac{\log\log
t}{\phi(n)} + O\Big( \frac{\log n}{\phi(n)} \Big).
\label{lose.the.minus.one}
\end{equation}
This readily implies that
\begin{equation}
\sum_{\substack{p\in\P n \\ p\le t}} \frac1{p-1}~ =~ \frac{\log\log
t}{\phi(n)} + O\Big( \frac{\log n}{\phi(n)} \Big)
\label{second.author}
\end{equation}
as well, since (noting that the smallest possible term in the sum is
$p=n+1$) the difference equals
\[
\sum_{\substack{p\in\P n \\ p\le t}} \frac1{(p-1)p} ~\le
~\sum_{i=1}^\infty \frac1{in(in+1)} \ll \frac1{n^2}.
\]

We occasionally use the Chebyshev upper bound
\begin{equation}
\sum_{p\le z} \log p~ \le~ \sum_{n\le z} \Lambda(n) ~\ll~ z,
\label{Chebyshev.bound}
\end{equation}
where $\Lambda(n)$ is the von Mangoldt function, as well as the weaker
versions \begin{equation}
\sum_{p\le z} \frac{\log p}p~ \ll ~\log z,\qquad \sum_{p\le z}
\frac{\log^2p}p~ \ll ~\log^2z
\label{weaker.Chebyshev.bound}
\end{equation}
and the tail estimates
\begin{equation}
\sum_{p>z} \frac{\log p}{p^2} ~\ll ~\frac1z,\qquad \sum_{p>z}  
\frac1{p^2}
\ll \frac1{z\log z},
\label{tail.estimates}
\end{equation}
each of which can be derived from the estimate \eqref{Chebyshev.bound}
by partial summation. We shall also need at one point a weak form of  
the asymptotic
formula of Mertens,
\begin{equation}
\sum_{p\le z} \frac{\log p}p~ =~ \log z+O(1).
\label{Mertens}
\end{equation}
For any polynomial $P(x)$, we also note the series estimate
\[
\sum_{a=0}^\infty \frac{P(a)}{m^a} ~\ll_P~ 1
\]
uniformly for $m\ge2$, valid since the series $\sum_{a=0}^\infty
P(a)z^a$ converges uniformly for $|z|\le\frac12$. The estimates
\begin{equation}
\sum_{a\in\N} \frac{P(a)}{m^a}~ \ll_P~ \frac1m, \qquad
\sum_{\substack{a\in\N \\ m^a>z}} \frac{P(a)}{m^a} ~\ll_P~ \frac1z,
\label{geometric.series}
\end{equation}
valid uniformly for any integer $m\ge2$, follow easily by factoring out
the first denominator occurring in each sum.

\section{Large primes dividing $\phi(\phi(n))$ and
$\lambda(\lambda(n))$}

\begin{proof}[Proof of Proposition \ref{same.large.prime.divisors.prop}]
If $q$ is any prime, then $q$ divides $\phi(\phi(n))$ if and only if at
least one of the following criteria holds:
\begin{itemize}
\item $q^3\mid n$,
\item there exists $p\in\P{q^2}$ with $p\mid n$,
\item there exists $p\in\P q$ with $p^2\mid n$,
\item there exist $r\in\P q$ and $p\in\P r$ with $p\mid n$,
\item $q^2\mid n$ and there exists $p\in\P q$ with $p\mid n$,
\item there exist distinct $p_1,p_2\in\P q$ with $p_1p_2\mid n$.
\end{itemize}
In the first four of these six cases, it is easily checked that
$q\mid\lambda(\lambda(n))$ as well. (This is not quite true for $q=2$,
but in this proof we shall only consider primes $q>y^2$.) Therefore we
can estimate the number of integers $n\le x$ for which $q$ divides
$\phi(\phi(n))$ but not $\lambda(\lambda(n))$ as follows:
\[
\sum_{\substack{n\le x \\ q\mid\phi(\phi(n)) \\
q\dnd\lambda(\lambda(n))}} 1
~\le~ \sum_{p\in\P q} \sum_{\substack{n\le x \\ q^2p\mid n}} 1 +
\sum_{p_1\in\P q} \sum_{\substack{p_2\in\P q \\ p_2\ne p_1}}
\sum_{\substack{n\le x \\ p_1p_2\mid n}} 1
~\le~ \sum_{p\in\P q} \frac x{q^2p} + \sum_{p_1\in\P q} \sum_{p_2\in\P
q} \frac x{p_1p_2}.
\]
Using three applications of the Brun--Titchmarsh inequality
\eqref{BT.prime.power}, we conclude that for any odd prime $q$,
\[
\sum_{\substack{n\le x \\ q\mid\phi(\phi(n)) \\
q\dnd\lambda(\lambda(n))}} 1~ \ll ~\frac{xy}{q^3} + \frac{xy^2}{q^2}  
~\ll
~\frac{xy^2}{q^2}.
\]
Consequently, by the tail estimate \eqref{tail.estimates} and the  
condition
$\psi(x)=o(\log y)$,
\begin{align*}
\sum_{q>y^2} \sum_{\substack{n\le x \\ q\mid\phi(\phi(n)) \\
q\dnd\lambda(\lambda(n))}} 1 ~\ll~ xy^2 \sum_{q>y^2} \frac1{q^2} ~\ll
~\frac{xy^2}{y^2\log y^2}~ <~ \frac x{\log y}~ \ll~ \frac x{\psi(x)}.
\end{align*}
Therefore for almost all $n\le x$, every prime $q>y^2$ dividing
$\phi(\phi(n))$ also divides $\lambda(\lambda(n))$, as asserted.
\end{proof}

\begin{lemma}
Given a real number $x\ge3$ and a prime $q>y^2$, define $S_q=S_q(x)$ to
be the set of all integers $n\le x$ for which at least one of the
following criteria holds:
\begin{itemize}
\item $q^2\mid n$,
\item there exists $p\in\P{q^2}$ with $p\mid n$,
\item there exist $r\in\P{q^2}$ and $p\in\P r$ with $p\mid n$,
\item there exist distinct $r_1,r_2,r_3\in\P q$ and $p\in\P{r_1r_2r_3}$
with $p\mid n$,
\item there exist distinct $r_1,r_2,r_3,r_4\in\P q$,
$p_1\in\P{r_1r_2}$, and $p_2\in\P{r_3r_4}$ with $p_1p_2\mid n$.
\end{itemize}
Then the cardinality of $S_q$ is $O(xy^2/q^2)$.
\label{less.squarefree.cases.lemma}
\end{lemma}

\begin{proof}
The number of integers up to $x$ for which any particular one of the
five criteria holds is easily shown to be $O(xy^2/q^2)$. For the sake
of conciseness, we show the details of this calculation only for the
last criterion, which is the most complicated. The number of integers  
$n$
up to $x$ for which there exist distinct $r_1,r_2,r_3,r_4\in\P q$,
$p_1\in\P{r_1r_2}$, and $p_2\in\P{r_3r_4}$ with $p_1p_2\mid n$ is at
most
\[
\sum_{r_1,r_2,r_3,r_4\in\P q} \sum_{\substack{p_1\in\P{r_1r_2} \\
p_2\in\P{r_3r_4}}} \sum_{\substack{n\le x \\ p_1p_2\mid n}} 1 ~\le
~\sum_{r_1,r_2,r_3,r_4\in\P q} \sum_{\substack{p_1\in\P{r_1r_2} \\
p_2\in\P{r_3r_4}}} \frac x{p_1p_2}.
\]
Using six applications of the Brun--Titchmarsh estimate
\eqref{BT.prime.power}, we have
\begin{equation*}
\sum_{r_1,r_2,r_3,r_4\in\P q} \sum_{\substack{p_1\in\P{r_1r_2} \\
p_2\in\P{r_3r_4}}} \frac x{p_1p_2} \ll \sum_{r_1,r_2,r_3,r_4\in\P q}
\frac{xy^2}{r_1r_2r_3r_4}~ \ll~ \frac{xy^6}{q^4} < \frac{xy^2}{q^2},
\end{equation*}
the last inequality being valid due to the hypothesis $q>y^2$.
\end{proof}

\begin{proof}[Proof of Proposition \ref{large.primes.in.phi.phi.prop}]
Define $S=S(x)$ to be the union of $S_q$ over all primes $q>y^2$, where
$S_q$ is defined as in the statement of Lemma
\ref{less.squarefree.cases.lemma}. Using $\#A$ to denote the
cardinality of a set $A$, Lemma \ref{less.squarefree.cases.lemma}
implies that
\[
\#S \le \sum_{q>y^2} \#S_q ~\ll~ \sum_{q>y^2} \frac{xy^2}{q^2}~ \ll
~\frac{xy^2}{y^2\log y^2} ~\ll~ \frac x{\psi(x)}
\]
by the tail estimate \eqref{tail.estimates} and the condition $\psi(x)
= o(\log y)$. Therefore to prove that the estimate
\eqref{large.primes.in.phi.phi.eq} holds for almost all integers $n\le
x$, it suffices to prove that it holds for almost all integers $n\le x$
that are not in the set $S$. This in turn is implied by the upper bound
\begin{equation}
\sum_{\substack{n\le x \\ n\notin S}} \sum_{\substack{q>y^2 \\
v_q(\phi(\phi(n)))\ge2}} v_q(\phi(\phi(n)))\log q ~\ll ~xy^2,
\label{not.counting.S.eq}
\end{equation}
which we proceed now to establish.

Fix a prime $q>y^2$ and an integer $a\ge2$ for the moment. In general,
there are many ways in which $q^a$ could divide $\phi(\phi(n))$,
depending on the power to which $q$ divides $n$ itself, the power to
which $q$ divides numbers of the form $p-1$ with $p\mid n$, and so
forth. However, for integers $n\notin S$, most of these various
possibilities are ruled out by one of the five criteria defining the
sets $S_q$. In fact, for $n\notin S$, there are only two ways for $q^a$
to divide $\phi(\phi(n))$:
\begin{itemize}
\item there are distinct $r_1,\dots,r_a\subset\P q$ and distinct
$p_1\in\P{r_1}$, \dots, $p_a\in\P{r_a}$ with $p_1\dots p_a|n$,
\item there are distinct $r_1,\dots,r_a\subset\P q$, distinct
$p_1\in\P{r_1}$, \dots,
$p_{a-2}\in\P{r_{a-2}}$, and $p\in\P{r_{a-1}r_a}$ with $p_1\dots
p_a|n$.
\end{itemize}
(We refer to the former case as the ``supersquarefree'' case.)

Still considering $q$ and $a$ fixed, the number of integers $n$ up to
$x$ satisfying each of these two conditions is at most
\[
\sum_{r_1,\dots,r_a\in\P q} \frac1{a!} \sum_{\substack{p_1\in\P{r_1} \\
\dots \\ p_a\in\P{r_a}}} \sum_{\substack{n\le x \\ p_1\dots p_a\mid n}}
1~\le~\sum_{r_1,\dots,r_a\in\P q} \frac1{a!}
\sum_{\substack{p_1\in\P{r_1} \\ \dots \\ p_a\in\P{r_a}}} \frac
x{p_1\dots p_a}
\]
and
\[
\sum_{r_1,\dots,r_a\in\P q} \frac1{2!(a-2)!}
\sum_{\substack{p_1\in\P{r_1} \\ \dots \\ p_{a-2}\in\P{r_{a-2}} \\
p\in\P{r_{a-1}r_a}}} \sum_{\substack{n\le x \\ p_1\dots p_{a-2}p\mid
n}} 1~\le~\sum_{r_1,\dots,r_a\in\P q} \frac1{(a-2)!}
\sum_{\substack{p_1\in\P{r_1} \\ \dots \\ p_{a-2}\in\P{r_{a-2}} \\
p\in\P{r_{a-1}r_a}}} \frac x{p_1\dots p_{a-2}p},
\]
respectively, the factors $1/a!$ and $1/2!(a-2)!$ coming from the
various possible permutations of the primes $r_i$. Letting $c\ge1$ be  
the constant implied in the
Brun--Titchmarsh inequality \eqref{BT.prime.power} as applied to moduli $n$
that are divisible by at most two distinct primes, we see that
\[
\sum_{r_1,\dots,r_a\in\P q} \frac1{a!} \sum_{\substack{p_1\in\P{r_1} \\
\dots \\ p_a\in\P{r_a}}} \frac x{p_1\dots p_a}~\le
~\sum_{r_1,\dots,r_a\in\P q} \frac1{a!} \frac{x(cy)^a}{r_1\dots r_a}~\le
~\frac{x(cy)^{2a}}{a!q^a}
\]
and
\[
\sum_{r_1,\dots,r_a\in\P q} \frac1{(a-2)!}
\sum_{\substack{p_1\in\P{r_1} \\ \dots \\ p_{a-2}\in\P{r_{a-2}} \\
p\in\P{r_{a-1}r_a}}} \frac x{p_1\dots p_{a-2}p}~\le
~\sum_{r_1,\dots,r_a\in\P q} \frac1{(a-2)!} \frac{x(cy)^{a-1}}{r_1\dots
r_a}~\le~\frac{x(cy)^{2a-1}}{(a-2)!q^a}.
\]
Therefore the number of integers $n\le x$ such that $n\notin S$ and
$q^a\mid\phi(\phi(n))$ is
\begin{equation}
\le~\frac{x(cy)^{2a}}{a!q^a} + \frac{x(cy)^{2a-1}}{(a-2)!q^a}~<
~\frac{c^{2a}xy^4}{(a-2)!q^2},
\label{two.case.bound}
\end{equation}
where we have used the assumption $q>y^2$.

We now establish the estimate \eqref{not.counting.S.eq}. Note that
\begin{align*}
\sum_{\substack{n\le x \\ n\notin S}}  
\sum_{\substack{q>y^2 \\
v_q(\phi(\phi(n)))\ge2}} v_q(\phi(\phi(n)))\log q
~&\le~
2\sum_{\substack{n\le x \\ n\notin S}}
\sum_{\substack{q>y^2 \\ v_q(\phi(\phi(n)))\ge2}} \big(
v_q(\phi(\phi(n))) -1\big) \log q \\
&=~2\sum_{q>y^2} \log q \sum_{a\ge2} \sum_{\substack{n\le x \\ n\notin S \\
q^a\mid\phi(\phi(n))}} 1.
\end{align*}
Therefore, using the bound \eqref{two.case.bound} for each pair $q$ and
$a$,
\begin{align*}
\sum_{\substack{n\le x \\ n\notin S}} \sum_{\substack{q>y^2 \\
v_q(\phi(\phi(n)))\ge2}} v_q(\phi(\phi(n)))\log q~&\le~2\sum_{q>y^2}
\log q \sum_{a\ge2} \frac{c^{2a}xy^4}{(a-2)!q^2} \\
&=~2c^4e^{c^2}xy^4 \sum_{q>y^2} \frac{\log  
q}{q^2}~\ll~\frac{xy^4}{y^2}~=~
xy^2
\end{align*}
by the tail estimate \eqref{tail.estimates}. This establishes the
estimate \eqref{not.counting.S.eq} and hence the proposition.
\end{proof}

\section{Small primes and the reduction to $h(n)$}

\begin{lemma}
For any prime power $q^a$, the number of positive integers $n\le x$ for
which $q^a$ divides $\lambda(\lambda(n))$ is $O(xy^2/q^a)$.
\label{prime.powers.in.lambda.lambda.lemma}
\end{lemma}

\begin{proof}
The prime power $q^a$ divides $\lambda(\lambda(n))$ only if at
least one of the following criteria holds:
\begin{itemize}
\item $q^{a+2}\mid n$,
\item there exists $p\in\P{q^a}$ with $p^2\mid n$,
\item there exists $p\in\P{q^{a+1}}$ with $p\mid n$,
\item there exist $r\in\P{q^a}$ and $p\in\P r$ with $p\mid n$.
\end{itemize}
Thus
\begin{align}
\sum_{\substack{n\le x \\ q^a\mid\lambda(\lambda(n))}} 1~ &\le
~\sum_{\substack{n\le x \\ q^{a+2}\mid n}} 1 + \sum_{p\in\P{q^a}}
\sum_{\substack{n\le x \\ p^2\mid n}} 1 + \sum_{p\in\P{q^{a+1}}}
\sum_{\substack{n\le x \\ p\mid n}} 1 + \sum_{r\in\P{q^a}} \sum_{p\in\P
r}
\sum_{\substack{n\le x \\ p\mid n}} 1 \notag \\
&\le~ \frac x{q^{a+2}} + \sum_{\substack{p\in\P{q^a} \\ p\le\sqrt x}}
\frac x{p^2} + \sum_{\substack{p\in\P{q^{a+1}} \\ p\le x}} \frac xp +
\sum_{r\in\P{q^a}} \sum_{\substack{p\in\P r\\p\le x}} \frac xp.
\label{four.criteria.split}
\end{align}
In the first of these three sums, it is sufficient to notice that any
$p\in\P{q^a}$ must exceed $q^a$, which leads to the estimate
\[
\sum_{\substack{p\in\P{q^a} \\ p\le\sqrt x}} \frac x{p^2}~ <~
\sum_{m>q^a} \frac x{m^2}~ <~ \frac x{q^a}.
\]
To bound the second and third sums, we invoke the Brun--Titchmarsh
estimate \eqref{BT.prime.power} a total of three times:
\begin{align*}
\sum_{\substack{p\in\P{q^{a+1}} \\ p\le x}} \frac xp ~&\ll~
\frac{xy}{q^{a+1}} \\
\sum_{r\in\P{q^a}} \sum_{\substack{p\in\P r\\p\le x}} \frac xp~ &\ll
~\sum_{\substack{r\in\P{q^a} \\ r\le x}} \frac{xy}r ~\ll
~\frac{xy^2}{q^a}.
\end{align*}
Using these three estimates, \eqref{four.criteria.split} gives
\[
\sum_{\substack{n\le x \\ q^a\mid\lambda(\lambda(n))}} 1 ~\ll~ \frac
x{q^{a+2}} + \frac x{q^a} + \frac{xy}{q^{a+1}} +\frac{xy^2}{q^a} ~\ll~
\frac{xy^2}{q^a},
\]
which establishes the lemma.
\end{proof}

\begin{proof}[Proof of Proposition
\ref{small.primes.in.lambda.lambda.prop}]
We have
\[
\sum_{q\le y^2} v_q(\lambda(\lambda(n))) \log q~ =~ \sum_{q\le y^2} \log
q \sum_{\substack{a\in\N \\ q^a\mid\lambda(\lambda(n))}} 1
~\le ~\sum_{q\le y^2} \log q \sum_{\substack{a\in\N \\ q^a\le y^2}} 1 +
\sum_{q\le y^2} \log q \sum_{\substack{a\in\N \\ q^a>y^2 \\
q^a\mid\lambda(\lambda(n))}} 1.
\]
Since the first sum is simply
\[
\sum_{q\le y^2} \log q \sum_{\substack{a\in\N \\ q^a\le y^2}} 1~ =~
\sum_{m\le y^2} \Lambda(m)~ \ll ~y^2
\]
by the Chebyshev estimate \eqref{Chebyshev.bound}, we have uniformly  
for $n\le x$,
\begin{equation}
\sum_{q\le y^2} v_q(\lambda(\lambda(n))) \log q ~\ll~ y^2 + \sum_{q\le
y^2} \log q \sum_{\substack{a\in\N \\ q^a>y^2 \\
q^a\mid\lambda(\lambda(n))}} 1.
\label{not.over.n.yet}
\end{equation}
To show that this quantity is usually small, we sum this last double
sum over $n$ and apply Lemma \ref{prime.powers.in.lambda.lambda.lemma},
yielding
\[
\sum_{n\le x} \sum_{q\le y^2} \log q \sum_{\substack{a\in\N \\ q^a>y^2
\\ q^a\mid\lambda(\lambda(n))}} 1~ =~ \sum_{q\le y^2} \log q
\sum_{\substack{a\in\N \\ q^a>y^2}} \sum_{\substack{n\le x \\
q^a\mid\lambda(\lambda(n))}} 1
~\ll ~\sum_{q\le y^2} \log q \sum_{\substack{a\in\N \\ q^a>y^2}}
\frac{xy^2}{q^a}.
\]
Using the geometric series sum \eqref{geometric.series} and the
Chebyshev estimate \eqref{Chebyshev.bound}, this becomes
\[
\sum_{n\le x} \sum_{q\le y^2} \log q \sum_{\substack{a\in\N \\ q^a>y^2
\\ q^a\mid\lambda(\lambda(n))}} 1 ~\ll ~\sum_{q\le y^2} \log q \cdot
\frac{xy^2}{y^2}~ \ll~ xy^2.
\]
Therefore if we sum both sides of \eqref{not.over.n.yet} over
$n$, we obtain
\begin{align*}
\sum_{n\le x} \sum_{q\le y^2} v_q(\lambda(\lambda(n))) \log q ~\ll  
~xy^2.
\end{align*}
This implies that for almost all $n\le x$, we have
\[
\sum_{q\le y^2} v_q(\lambda(\lambda(n))) \log q ~\ll~ y^2\psi(x),
\]
as desired.
\end{proof}

\begin{proof}[Proof of Proposition \ref{h.good.approximation.prop}]
Fix a prime $q$ for the moment. For any positive integer $m$, the
usual formula for $\phi(m)$ readily implies
\[
v_q(\phi(m))~ =~ \max\{0,v_q(m)-1\} + \sum_{p\mid m} v_q(p-1),
\]
which we use in the form
\[
\sum_{p\mid m} v_q(p-1) ~\le~ v_q(\phi(m)) ~\le~ \sum_{p\mid m}  
v_q(p-1) +
v_q(m).
\]
Using these inequalities twice, first with $m=\phi(n)$ and then with
$m=n$, we see that
\begin{align}
\sum_{p\mid\phi(n)} v_q(p-1)~ \le~ v_q(\phi(\phi(n)))~ &\le
~\sum_{p\mid\phi(n)} v_q(p-1) + v_q(\phi(n)) \notag \\
&\le~ \sum_{p\mid\phi(n)} v_q(p-1) + \sum_{p\mid n} v_q(p-1) + v_q(n).
\label{two.phi.inequalities}
\end{align}
Now a prime $r$ divides $\phi(n)$ if and only if either $r^2\mid n$ or
there exists a prime $p\mid n$ such that $r\mid p-1$. Therefore
\[
\sum_{p\mid n} \sum_{r\mid p-1} v_q(r-1) ~\le~ \sum_{r\mid\phi(n)}
v_q(r-1)~ \le~ \sum_{p\mid n} \sum_{r\mid p-1} v_q(r-1) + \sum_{r\colon
r^2\mid n} v_q(r-1),
\]
the latter inequality accounting for the possibility that both criteria
hold for some prime~$r$. When we combine these inequalities with those
in equation \eqref{two.phi.inequalities} and subtract the double sum
over $p$ and $r$ throughout, we obtain
\begin{align*}
0~ \le~ v_q(\phi(\phi(n))) - \sum_{p\mid n} \sum_{r\mid p-1} v_q(r-1)
~&\le~ \sum_{r\colon r^2\mid n} v_q(r-1) + \sum_{p\mid n} v_q(p-1) +
v_q(n) \\
&\le~ 2\sum_{p\mid n} v_q(p-1) + v_q(n).
\end{align*}
Now we multiply through by $\log q$ and sum over all primes $q\le y^2$
to conclude that for any positive integer $n$,
\begin{equation*}
0~ \le~ \sum_{q\le y^2} v_q(\phi(\phi(n))) \log q - h(n) ~\le~ 2  
\sum_{q\le
y^2} \sum_{p\mid n} v_q(p-1) \log q + \sum_{q\le y^2} v_q(n) \log q.
\end{equation*}

It remains to show that the right-hand side of this last inequality is
$O(y\log y\cdot\psi(x))$ for almost all $n\le x$, which we accomplish
by establishing the estimate
\begin{equation}
\sum_{n\le x} \sum_{q\le y^2} \sum_{p\mid n} v_q(p-1) \log q +
\sum_{n\le x} \sum_{q\le y^2} v_q(n) \log q ~\ll~ xy\log y.
\label{hn.prop.to.show}
\end{equation}
We may rewrite the first term on the left-hand side as
\begin{align*}
\sum_{n\le x} \sum_{q\le y^2} \sum_{p\mid n} v_q(p-1) \log q~ &=
~\sum_{n\le x} \sum_{q\le y^2} \sum_{p\mid n} \sum_{\substack{a\in\N \\
q^a\mid p-1}} \log q \\
&=~ \sum_{q\le y^2} \log q \sum_{a\in\N} \sum_{p\in\P{q^a}}
\sum_{\substack{n\le x \\ p\mid n}} 1 ~\le~ \sum_{q\le y^2} \log q
\sum_{a\in\N} \sum_{p\in\P{q^a}} \frac xp.
\end{align*}
Using the Brun--Titchmarsh inequality \eqref{BT.prime.power} and the
geometric series estimate \eqref{geometric.series}, we obtain
\[
\sum_{n\le x} \sum_{q\le y^2} \sum_{p\mid n} v_q(p-1) \log q ~\ll~ x
\sum_{q\le y^2} \log q \sum_{a\in\N} \frac y{q^a}~ \ll~ xy \sum_{q\le
y^2} \frac{\log q}q~ \ll~ xy\log y^2.
\]
The second term on the left-hand side of \eqref{hn.prop.to.show} is
even simpler: we have
\[
\sum_{n\le x} \sum_{q\le y^2} v_q(n) \log q~ =~ \sum_{q\le y^2} \log q
\sum_{a\in\N} \sum_{\substack{n\le x \\ q^a\mid n}} 1 ~\le~ \sum_{q\le
y^2} \log q \sum_{a\in\N} \frac x{q^a},
\]
and using the geometric series bound \eqref{geometric.series} and the
weak Chebyshev estimate \eqref{weaker.Chebyshev.bound} yields
\[
\sum_{n\le x} \sum_{q\le y^2} v_q(n) \log q~ \ll~ x \sum_{q\le y^2}
\frac{\log q}q ~\ll~ x\log y^2.
\]
The last two estimates therefore establish \eqref{hn.prop.to.show} and
hence the proposition.
\end{proof}

\section{The normal order of $h(n)$}\label{h.normal.order.section}

Recall the definition \eqref{h.definition}:
$h(n)~ =~ \sum_{p\mid n} \sum_{r\mid p-1} \sum_{q\le y^2} v_q(r-1) \log  
q$.
We now calculate the normal order of the additive function $h(n)$ via
the Tur\'an--Kubilius inequality (see \cite{K}, Lemma 3.1). If we
define
\[
M_1(x)~=~\sum_{p\le x}\frac{h(p)}p,\qquad M_2(x)~=~\sum_{p\le
x}\frac{h(p)^2}p,
\]
then the Tur\'an-Kubilius inequality asserts that
\begin{equation}
\sum_{n\le x}(h(n)-M_1(x))^2~\ll ~xM_2(x).
\label{TK}
\end{equation}
\begin{prop}
We have $M_1(x) = y^2\log y + O(y^2)$ for all $x>e^{e^e}$.
\label{M1.prop}
\end{prop}

\begin{prop}
We have $M_2(x) \ll y^3\log^2y$ for all $x>e^{e^e}$.
\label{M2.prop}
\end{prop}

\begin{proof}[Proof of Proposition \ref{h.normal.order.prop}]
Let $N$ denote the number of $n\le x$ for which
$|h(n)-M_1(x)|>y^2$.  The contribution of such $n$ to the sum
in \eqref{TK} is at least $y^4N$.  Thus, Proposition \ref{M2.prop}
implies that $N\ll x(\log y)^2/y$.  Hence, Proposition \ref{M1.prop}
implies that $h(n)=y^2\log y+O(y^2)$ for all $n\le x$ but for a
set of size $O(x(\log y)^2)/y)$.  This proves Proposition  
\ref{h.normal.order.prop}.
\end{proof}

To calculate $M_1(x)$ and $M_2(x)$ we shall first calculate $\sum_{p\le
t} h(p)$ and $\sum_{p\le t} h(p)^2$ and then account for the weights
$1/p$ using partial summation. We begin the evaluation of $\sum_{p\le
t} h(p)$ with a lemma.

\begin{lemma}
Let $b$ be a positive integer and $t>e^e$ a real number.
\begin{enumerate}
\item[(a)] If $b>t^{1/4}$ then
\[
\sum_{r\in\P b} \pi(t;r,1)~ \ll~ \frac{t\log t}b.
\]
\item[(b)] If $b\le t^{1/4}$ then
\[
\sum_{\substack{r\in\P b \\ r>t^{1/3}}} \pi(t;r,1) ~\ll~ \frac
{bt}{\phi(b)^2\log t}.
\]
and
\[
\sum_{r\in\P b} \pi(t;r,1) ~\ll~ \frac{t\log\log t}{\phi(b)\log t}
\]
\end{enumerate}
\label{tricky.to.get.right.lemma}
\end{lemma}

\noindent {\em Remark.} The exponents $\frac14$ and $\frac13$ are
rather arbitrary and chosen only for simplicity; any two exponents
$0<\alpha<\beta<\frac12$ would do equally well.

\begin{proof}
Notice that in all three sums, the only contributing terms are those
with $r>b$ and $r<t$. If $b>t^{1/4}$, then the trivial bound
$\pi(t;r,1)\le t/r$ gives
\begin{equation*}
\sum_{r\in\P{b}}\pi(t;r,1) ~\le~ \sum_{\substack{r\in\P{b} \\
t^{1/4}<r\le t}} \frac tr ~\le
~\sum_{\substack{m\equiv1\mod{b}\\ t^{1/4}<m\le t}}\frac tm
~\ll~\frac{t\log t}{b},
\end{equation*}
proving part (a) of the lemma.

We now assume $b\le t^{1/4}$. We have
\begin{align*}
\sum_{\substack{r\in\P b \\ r>t^{1/3}}} \pi(t;r,1)~ &=~ \#\{ (m,r)\colon
r\equiv1\mod{b},\, r>t^{1/3},\, mr+1\le t,\, \text{$mr+1$ and $r$ both
prime} \} \\
&\le~ \sum_{m<t^{2/3}} \#\{ r<\tfrac tm\colon r\equiv1\mod{b},\,
\text{$mr+1$ and $r$ both prime} \} \\
&\ll~ \sum_{m<t^{2/3}} \frac bt{\phi(mb)\phi(b)\log^2 \frac t{mb}}
\end{align*}
by Brun's sieve method (see \cite[Corollary 2.4.1]{HR}).
We have $\frac t{mb} \ge t^{1/12}$ and so $\log
\frac t{mb} \gg \log t$. We also have $\phi(mb) \ge \phi(m)\phi(b)$ and
the standard estimate
\begin{equation}
\sum_{m\le z} \frac1{\phi(m)} ~\ll~ \log z.
\label{phi.reciprocal.sum}
\end{equation}
Therefore
\[
\sum_{\substack{r\in\P b \\ r>t^{1/3}}} \pi(t;r,1) ~\ll~  
\sum_{m<t^{2/3}}
\frac {bt}{\phi(m)\phi(b)^2\log^2 t}  ~\ll~ \frac{bt\log
t^{2/3}}{\phi(b)^2\log^2t} ~\le~ \frac {bt}{\phi(b)^2\log t},
\]
establishing the first estimate in part (b). Finally, by the
Brun--Titchmarsh inequalities \eqref{real.BT} and
\eqref{BT.prime.power},
\begin{equation*}
\sum_{\substack{r\in\P b \\ r\le t^{1/3}}} \pi(t;r,1) ~\ll~
\sum_{\substack{r\in\P b \\ r\le t^{1/3}}} \frac t{\phi(r)\log\frac tr}
~\ll~ \sum_{\substack{r\in\P b \\ r\le t^{1/3}}} \frac t{r\log t}~ \ll~
\frac {t\log\log t}{\phi(b)\log t}.
\end{equation*}
Combining this estimate with the first half of part (b) and the standard
estimate $b/\phi(b)\ll\log\log b$ establishes the second half.
\end{proof}

\begin{lemma}
For all real numbers $x>e^{e^e}$ and $t>e^e$, we have
\begin{equation*}
\sum_{p\le t}h(p)~ = ~\frac{2t\log\log t\log y}{\log t} + O\Big(
\frac{t\log\log t}{\log t} + \frac{t\log^2y}{\log t} + t^{3/4}\log
t\cdot y^2 \Big).
\end{equation*}
\label{hp.sum.lemma}
\end{lemma}

\noindent{\em Remark.} In particular, we have $\sum_{p\le x}h(p) \ll
x\log\log x\log y/\log x = xy\log y/\log x$.

\begin{proof}
We may rewrite
\begin{multline}
\sum_{p\le t} h(p)~ =~ \sum_{p\le t} \sum_{r\mid p-1} \sum_{q\le y^2}
v_q(r-1) \log q~ =~ \sum_{p\le t} \sum_{r\mid p-1} \sum_{q\le y^2}
\sum_{\substack{a\in\N \\ q^a \mid r-1}} \log q \\
~=~ \sum_{q\le y^2} \log q \sum_{a\in\N} \sum_{r\colon q^a \mid r-1}
\sum_{\substack{p\le t \\ r\mid p-1}} 1 ~=~ \sum_{q\le y^2} \log q
\sum_{a\in\N} \sum_{r\in\P{q^a}} \pi(t;r,1).
\label{hp.sum.before.split}
\end{multline}
The main contribution to this triple sum comes from the terms with
$q^a\le t^{1/4}$ and $r\le t^{1/3}$.
In fact, using Lemma \ref{tricky.to.get.right.lemma}(a) we can
bound the contribution from the terms with $q^a$ large by
\begin{align*}
\sum_{q\le y^2} \log q \sum_{\substack{a\in\N \\ q^a>t^{1/4}}}
\sum_{r\in\P{q^a}} \pi(t;r,1) ~&\le~ \sum_{q\le
y^2} \log q \sum_{\substack{a\in\N \\ q^a>t^{1/4}}} \sum_{r\in\P{q^a}}
\pi(t;r,1)
~\ll~ \sum_{q\le y^2} \log q \sum_{\substack{a\in\N \\ q^a>t^{1/4}}}
\frac{t\log t}{q^a} \\
&\ll~ t\log t \sum_{q\le y^2} \frac{\log q}{t^{1/4}}~ \ll~ t^{3/4}\log
t\cdot y^2,
\end{align*}
where the last two estimates are due to the geometric series bound
\eqref{geometric.series} and the Chebyshev bound
\eqref{Chebyshev.bound}. Similarly, using the first half of Lemma
\ref{tricky.to.get.right.lemma}(b) we can bound the contribution from
the terms with $q^a$ small and $r$ large by
\[
\sum_{q\le y^2} \log q \sum_{\substack{a\in\N \\ q^a\le t^{1/4}}}
\sum_{\substack{r\in\P{q^a} \\ r>t^{1/3}}} \pi(t;r,1)~\ll~ \sum_{q\le
y^2} \log q \sum_{\substack{a\in\N \\ q^a\le t^{1/4}}} \frac
t{q^a\log t}
~\ll~\frac t{\log t}\sum_{q\le y^2}\frac{\log q}q~\ll~\frac{t\log  
y}{\log t},
\]
where again the last two estimates are due to the geometric series bound
\eqref{geometric.series} and the weak Chebyshev bound
\eqref{weaker.Chebyshev.bound}. In light of these two estimates,
equation \eqref{hp.sum.before.split} becomes
\begin{equation}
\sum_{p\le t} h(p) ~=~ \sum_{q\le y^2} \log q \sum_{\substack{a\in\N\\  
q^a\le t^{1/4}}}
\sum_{\substack{r\in\P{q^a} \\ r\le t^{1/3}}} \pi(t;r,1) + O \Big(
t^{3/4}\log t\cdot y^2 + \frac{t\log y}{\log t} \Big).
\label{hp.sum.after.split}
\end{equation}

Define $E(t;r,1) = \pi(t;r,1) - \li(t)/(r-1)$. We have
\begin{multline}
\sum_{q\le y^2} \log q \sum_{\substack{a\in\N\\ q^a\le t^{1/4}}}  
\sum_{\substack{r\in\P{q^a} \\
r\le t^{1/3}}} \pi(t;r,1)~=~\sum_{q\le y^2} \log q  
\sum_{\substack{a\in\N\\ q^a\le t^{1/4}}}
\sum_{\substack{r\in\P{q^a} \\ r\le
t^{1/3}}} \Big( \frac{\li(t)}{r-1} + E(t;r,1) \Big) \\
=~\sum_{q\le y^2} \log q \sum_{\substack{a\in\N\\ q^a\le t^{1/4}}}
\sum_{\substack{r\in\P{q^a} \\ r\le t^{1/3}}} \frac{\li(t)}{r-1}
+ O\bigg( \sum_{q\le y^2} \log q \sum_{\substack{a\in\N\\ q^a\le  
t^{1/4}}}
\sum_{\substack{r\in\P{q^a} \\ r\le t^{1/3}}} |E(t;r,1)| \bigg).
\label{set.up.BV.application}
\end{multline}
Let $\Omega(m)$ denote the number of divisors of $m$ that are primes or  
prime powers.
Using the estimate $\Omega(m) \ll \log m$, we quickly dispose of
\begin{align*}
\sum_{q\le y^2} \log q \sum_{\substack{a\in\N\\ q^a\le t^{1/4}}}  
\sum_{\substack{r\in\P{q^a} \\
r\le t^{1/3}}} |E(t;r,1)|~&=~\log y \sum_{r\le t^{1/3}} |E(t;r,1)|
\sum_{q\le y^2} \sum_{\substack{a\in\N \\ q^a\mid r-1}} 1 \\
&\le~\log y \sum_{r\le t^{1/3}} |E(t;r,1)|\, \Omega(r-1) \\
&\ll~\log y\log t\sum_{r\le t^{1/3}}|E(t;r,1)|~\ll~\frac{t\log y}{\log
t}
\end{align*}
by the Bombieri--Vinogradov theorem (we could equally well put any
power of $\log t$ in the denominator of the final expression if we
needed). Inserting this estimate into equation
\eqref{set.up.BV.application}, we see that equation
\eqref{hp.sum.after.split} becomes
\begin{equation}
\sum_{p\le t} h(p)~=~\li(t) \sum_{q\le y^2} \log q \sum_{a\in\N}
\sum_{\substack{r\in\P{q^a} \\ r\le t^{1/3}}} \frac1{r-1} + O \Big(
t^{3/4}\log t\cdot y^2 + \frac{t\log y}{\log t} \Big).
\label{hp.sum.after.BV}
\end{equation}

We have by equation \eqref{second.author}
\begin{align*}
\sum_{q\le y^2} \log q&\sum_{a\in\N} \sum_{\substack{r\in\P{q^a} \\
r\le t^{1/3}}} \frac1{r-1}~=~\sum_{q\le y^2} \log q \sum_{a\in\N}
\Big( \frac{\log\log
t^{1/3}}{\phi(q^a)} + O\Big( \frac{\log q^a}{q^a} \Big) \Big) \\
&=~(\log\log t + O(1)) \sum_{q\le y^2} \log q \sum_{a\in\N} \Big(
\frac1{q^a} + O\Big( \frac1{q^{a+1}} \Big) \Big)
+ O\bigg( \sum_{q\le y^2} \log^2 q \sum_{a\in\N} \frac{a}{q^a}
\bigg) \\
&=~(\log\log t + O(1)) \sum_{q\le y^2} \Big( \frac{\log q}{q} + O\Big(
\frac{\log q}{q^2} \Big) \Big) + O\bigg( \sum_{q\le y^2} \frac{\log^2
q}q \bigg),
\end{align*}
using the geometric series estimate \eqref{geometric.series}. Using the
Mertens formula \eqref{Mertens} to evaluate the main term and the weak
Chebyshev estimates \eqref{weaker.Chebyshev.bound} to bound the error
terms, we see that
\[
\sum_{q\le y^2} \log q \sum_{a\in\N} \sum_{\substack{r\in\P{q^a} \\
r\le t^{1/3}}} \frac1{r-1}~=~\log\log t\log y^2 + O(\log y + \log\log t
+ \log^2y).
\]
We conclude from equation \eqref{hp.sum.after.BV} and the fact that
$\li(t) = t/\log t + O( t/\log^2t)$ that
\begin{align*}
\sum_{p\le t} h(p)~&=~\li(t) \big( \log\log t\log y^2 + O(\log y +
\log\log t + \log^2y) \big) \\
&\qquad{}+ O \Big( t^{3/4}\log t\cdot y^2 + \frac{t\log y}{\log t}
\Big) \\
&=~\frac{2t\log\log t\log y}{\log t} + O\Big( \frac{t\log\log t}{\log
t} + \frac{t\log^2y}{\log t} + t^{3/4}\log t\cdot y^2 \Big),
\end{align*}
as asserted.
\end{proof}

\begin{proof}[Proof of Proposition \ref{M1.prop}]
In an explicit example of the technique of partial summation, we write
\begin{align*}
M_1(x)~=~\sum_{p\le x} \frac{h(p)}p~&=~\sum_{p\le e^e} \frac{h(p)}p +
\sum_{e^e<p\le x} h(p) \bigg( \frac1x + \int_p^x \frac{dt}{t^2} \bigg)
\\
&=~O(1) + \frac1x \sum_{e^e<p\le x} h(p) + \int_{e^e}^x \frac{dt}{t^2}
\sum_{e^e<p\le t} h(p).
\end{align*}
The quantity $\sum_{p\le t} h(p)$ has been evaluated asymptotically in
Lemma \ref{hp.sum.lemma}, and the quantity $\sum_{e^e<p\le t} h(p)$
differs by only $O(1)$. Therefore we may use Lemma \ref{hp.sum.lemma}
and the remark following its statement to write
\begin{align*}
M_1(x)~&=~O(1) + \frac1x O\Big( \frac{xy\log y}{\log x} \Big) \\
&\qquad{}+ \int_{e^e}^x \frac{dt}{t^2} \Big( \frac{2t\log\log t\log
y}{\log t} + O\Big( \frac{t\log\log t}{\log t} + \frac{t\log^2y}{\log
t} + t^{3/4}\log t\cdot y^2 \Big) \Big) \\
&=~O\Big( \frac{y\log y}{\log x} \Big) +\log y \int_{e^e}^x
\frac{2\log\log t}{t\log t} \,dt \\
&\qquad{}+ O\bigg(  \int_{e^e}^x \frac{\log\log t}{t\log t} \,dt +
\log^2y \int_{e^e}^x \frac{dt}{t\log t} +  y^2 \int_{e^e}^x
\frac{dt}{t^{5/4}} \bigg).
\end{align*}
Each of these integrals can be explicitly evaluated, resulting in the
asymptotic formula
\begin{align*}
M_1(x)~&=~\log y \big( (\log\log x)^2-1 \big)+ O\Big( \frac{y\log
y}{\log x} + (\log\log x)^2 + \log^2y\cdot \log\log x + y^2 \Big) \\
&=~y^2\log y + O(y^2),
\end{align*}
as claimed.
\end{proof}

Now we turn our attention to $M_2(x)$, beginning with some preliminary
lemmas.

\begin{lemma}
For all real numbers $x>e^{e^e}$ and $t>e^e$, we have
\[
\sum_{q_1,q_2\le y^2} \log q_1 \log q_2 \sum_{a_1,a_2\in\N}
\sum_{r\in\P{q_1^{a_1}} \cap \P{q_2^{a_2}}} \sum_{\substack{p\le t \\
p\equiv1\mod{r}}} 1~\ll~t^{7/8}\log t\cdot y^2\log y + \frac{t\log\log
t\cdot \log^2y}{\log t}.
\]
\label{hp2.sum.error.lemma}
\end{lemma}

\begin{proof}
Since the exact form of $\P{q_1^{a_1}} \cap \P{q_2^{a_2}}$ depends on
whether or not $q_1=q_2$, we split the expression in question into two
separate sums:
\begin{align}
\sum_{q_1,q_2\le y^2}&\log q_1 \log q_2\sum_{a_1,a_2\in\N}
\sum_{r\in\P{q_1^{a_1}} \cap \P{q_2^{a_2}}} \sum_{\substack{p\le t \\
p\equiv1\mod{r}}} 1
\label{qs.equal.or.not} \\
&=~\sum_{q\le y^2} \log^2q \sum_{a_1,a_2\in\N}
\sum_{r\in\P{q^{\max\{a_1,a_2\}}}} \pi(t;r,1)
+ \sum_{\substack{q_1,q_2\le y^2 \\ q_1\ne q_2}} \log q_1 \log
q_2 \sum_{a_1,a_2\in\N} \sum_{r\in\P{q_1^{a_1}q_2^{a_2}}} \pi(t;r,1).
\notag
\end{align}
Noting that there are exactly $2a-1$ ordered pairs $(a_1,a_2)$ for
which $\max\{a_1,a_2\}=a$, we have
\begin{align*}
\sum_{q\le y^2} \log^2q&\sum_{a_1,a_2\in\N}
\sum_{r\in\P{q^{\max\{a_1,a_2\}}}}
\pi(t;r,1)~=~\sum_{q\le y^2} \log^2q \sum_{a\in\N} (2a-1)
\sum_{r\in\P{q^a}} \pi(t;r,1) \\
&\ll~\sum_{q\le y^2} \log^2q \sum_{\substack{a\in\N \\
q^a>t^{1/4}}}\frac{at\log t}{q^a}
+ \sum_{q\le y^2} \log^2q \sum_{\substack{a\in\N \\ q^a\le
t^{1/4}}}\frac{at\log\log t}{q^a\log t}
\end{align*}
by Lemma \ref{tricky.to.get.right.lemma}. Since
\[
\sum_{q\le y^2} \log^2q \sum_{\substack{a\in\N \\
q^a>t^{1/4}}}\frac{at\log t}{q^a}~\ll~t\log t \log y^2\sum_{q\le y^2}
\frac{\log q}{t^{1/4}}~\ll~t^{3/4}\log t\cdot y^2\log y
\]
by the Chebyshev bound \eqref{Chebyshev.bound}, and
\[
\sum_{q\le y^2} \log^2q \sum_{\substack{a\in\N \\ q^a\le
t^{1/4}}}\frac{at\log\log t}{q^a\log t}~\ll~\frac{t\log\log
t}{\log t} \sum_{q\le y^2} \frac{\log^2q}q~\ll~\frac{t\log\log
t\cdot\log^2y}{\log t}
\]
by \eqref{Chebyshev.bound} and its weaker version
\eqref{weaker.Chebyshev.bound}, the first term on the right-hand side
of equation \eqref{qs.equal.or.not} is bounded by the estimate asserted
in the statement of the lemma.

It remains to satisfactorily bound the second term on the right-hand
side of equation \eqref{qs.equal.or.not}. Again dividing the sum so
that Lemma \ref{tricky.to.get.right.lemma} can be applied, we have
\begin{align*}
\sum_{\substack{q_1,q_2\le y^2 \\ q_1\ne q_2}} \log q_1 \log q_2
\sum_{a_1,a_2\in\N} \sum_{r\in\P{q_1^{a_1}q_2^{a_2}}} & \pi(t;r,1)~\ll~
\sum_{q_1,q_2\le y^2} \log q_1 \log q_2 \sum_{\substack{a_1,a_2\in\N \\
q_1^{a_1}q_2^{a_2} > t^{1/4}}} \frac{t\log t}{q_1^{a_1}q_2^{a_2}} \\
&\qquad{}+ \sum_{q_1,q_2\le y^2} \log q_1 \log q_2
\sum_{\substack{a_1,a_2\in\N \\ q_1^{a_1}q_2^{a_2}\le t^{1/4}}}
\frac{t\log\log t}{q_1^{a_1}q_2^{a_2}\log t}.
\end{align*}
In the first of these two terms, at least one of the $q_i^{a_i}$ must
exceed $t^{1/8}$, and so using the estimates \eqref{geometric.series},
\eqref{Chebyshev.bound}, and \eqref{weaker.Chebyshev.bound} we see that
\begin{align*}
\sum_{q_1,q_2\le y^2} \log q_1 \log q_2 \sum_{\substack{a_1,a_2\in\N \\
q_1^{a_1}q_2^{a_2} > t^{1/4}}} \frac{t\log t}{q_1^{a_1}q_2^{a_2}}~&\le~
2t\log t \sum_{q_1\le y^2} \log q_1 \sum_{\substack{a_1\in\N \\  
q_1^{a_1} >
t^{1/8}}} \frac1{q_1^{a_1}} \sum_{q_2\le y^2} \log q_2 \sum_{a_2\in\N}
\frac1{q_2^{a_2}} \\
&\ll~t\log t \sum_{q_1\le y^2} \frac{\log q_1}{t^{1/8}} \sum_{q_2\le
y^2} \frac{\log q_2}{q_2} \\
&\ll~t^{7/8}\log t\cdot y^2\log y.
\end{align*}
In the second, we simply ignore the restriction $q_1^{a_1}q_2^{a_2}\le
t^{1/4}$ and use the estimates \eqref{geometric.series} and
\eqref{weaker.Chebyshev.bound}, obtaining
\begin{align*}
\sum_{q_1,q_2\le y^2} \log q_1
\log q_2 \sum_{a_1,a_2\in\N} \frac{t\log\log t}{q_1^{a_1}q_2^{a_2}\log
t}
&=~\frac{t\log\log t}{\log t} \bigg( \sum_{q\le y^2} \log q
\sum_{a\in\N} \frac1{q^a} \bigg)^2 \\
&\ll~\frac{t\log\log t}{\log t} \bigg( \sum_{q\le y^2} \frac{\log q}q
\bigg)^2 \\
&\ll~\frac{t\log\log t\cdot\log^2y}{\log t}.
\end{align*}
This concludes the proof of the lemma.
\end{proof}

The following lemma is similar in spirit to Lemma  
\ref{tricky.to.get.right.lemma} but is a bit more complicated to state  
and prove.

\begin{lemma}
Let $b_1$ and $b_2$ be positive integers and $t>e^e$ a real number.
\begin{enumerate}
\item[(a)] If $b_1>t^{1/8}$ or $b_2>t^{1/8}$ then
\[
\sum_{r_1\in\P{b_1}} \sum_{r_2\in\P{b_2}} \pi(t;r_1r_2,1)~\ll~
\frac{t\log^2t}{b_1b_2}.
\]
\item[(b)] If neither $b_1$ nor $b_2$ exceeds $t^{1/8}$ then
\[
\sum_{r_1\in\P{b_1}} \sum_{\substack{r_2\in\P{b_2} \\ r_1r_2>t^{1/3}}}
\pi(t;r_1r_2,1)~\ll~\frac{b_2t\log\log t}{\phi(b_1)\phi(b_2)^2\log t}
\]
and
\[
\sum_{r_1\in\P{b_1}} \sum_{r_2\in\P{b_2}} \pi(t;r_1r_2,1)~\ll~
\frac{t(\log\log t)^2}{\phi(b_1)\phi(b_2)\log t}.
\]
\end{enumerate}
\label{double.tricky.lemma}
\end{lemma}

\noindent {\em Remark.} Again, the values $1/8$ and $1/3$ for the  
exponents
are rather arbitrary.

\begin{proof}
The bound in part (a) follows from the trivial estimate
$\pi(t;r_1r_2,1) \ll t/r_1r_2$, just as in the proof of Lemma
\ref{tricky.to.get.right.lemma}(a). For the first estimate in part (b),
we my assume that $r_1\le r_2$ by symmetry. We use Brun's method again:
\begin{align*}
\sum_{r_1\in\P{b_1}} \sum_{\substack{r_2\in\P{b_2} \\ r_1\le r_2 \\
r_1r_2>t^{1/3}}} &\pi(t;r_1r_2,1) \\
&=~\#\{(m,r_1,r_2)\colon r_1\equiv1\mod{b_1},\, r_2\equiv1\mod{b_2},\,
r_1\le r_2,\, r_1r_2>t^{1/3}, \\
&\qquad\quad mr_1r_2+1\le t,\, \text{and $r_1$, $r_2$, and $mr_1r_2+1$
are all prime} \} \\
&\le~\sum_{m<t^{2/3}} \sum_{\substack{r_1<\sqrt{t/m}\\ r_1\in\P{b_1} }}
\sum_{\substack{r_2<t/mr_1\\ r_2\in\P{b_2}\\ mr_1r_2+1~{\rm prime}}}1\\
%
%
&\ll~\sum_{m<t^{2/3}} \sum_{\substack{r_1<\sqrt{t/m}\\ r_1\in\P{b_1} }}
\frac{mr_1b_2}{\phi(b_2)\phi(mr_1b_2)}\cdot\frac{t/mr_1}{\log^2(t/ 
mr_1b_2)}.
\end{align*}
Notice that $t/mr_1b_2 > (\sqrt{t/m})/b_2 > t^{1/6}/t^{1/8} =
t^{1/24}$, and so
\begin{align*}
\sum_{r_1\in\P{b_1}} \sum_{\substack{r_2\in\P{b_2}\\ r_1\le r_2\\
r_1r_2>t^{1/3}}} \pi(t;r_1r_2,1)
&\ll~\frac t{\log^2t}\sum_{m<t^{2/3}}\sum_{\substack{r_1<\sqrt{t/m}\\  
r_1\in\P{b_1}}}
\frac{b_2}{\phi(b_2)^2\phi(m)\phi(r_1)}\\
&\ll~\frac{b_2t\log\log  
t}{\phi(b_1)\phi(b_2)^2\log^2t}\sum_{m<t^{2/3}}\frac1{\phi(m)}\\
&\ll~\frac{b_2t\log\log t}{\phi(b_1)\phi(b_2)^2\log t}.
%
\end{align*}
by the estimates \eqref{primesum} and \eqref{phi.reciprocal.sum} as
desired. The second estimate of part (b) is a consequence of the first  
estimate and
\[
\sum_{r_1\in\P{b_1}} \sum_{\substack{r_2\in\P{b_2} \\ r_1r_2\le
t^{1/3}}} \pi(t;r_1r_2,1)~\ll~\frac{t(\log\log
t)^2}{\phi(b_1)\phi(b_2)\log t},
\]
which follows from the Brun--Titchmarsh inequality just as in the proof
of Lemma \ref{tricky.to.get.right.lemma}(b).
\end{proof}

\begin{proof}[Proof of Proposition \ref{M2.prop}]
We may rewrite
\begin{align*}
\sum_{p\le t} h(p)^2~&=~\sum_{p\le t} \bigg( \sum_{r\mid
p-1} \sum_{q\le y^2} \sum_{\substack{a\in\N \\ q^a\mid r-1}} \log q
\bigg)^2 \\
&=~\sum_{q_1,q_2\le y^2} \log q_1 \log q_2 \sum_{a_1,a_2\in\N}
\sum_{\substack{r_1\in\P{q_1^{a_1}} \\ r_2\in\P{q_2^{a_2}}}}
\sum_{\substack{p\le t \\ p\equiv1\mod{r_1} \\ p\equiv1\mod{r_2}}} 1 \\
&=~\sum_{q_1,q_2\le y^2} \log q_1 \log q_2 \sum_{a_1,a_2\in\N}
\sum_{\substack{r_1\in\P{q_1^{a_1}} \\ r_2\in\P{q_2^{a_2}} \\ r_1\ne
r_2}} \sum_{\substack{p\le t \\ p\equiv1\mod{r_1} \\
p\equiv1\mod{r_2}}} 1 \\
&\qquad{}+ O\Big( t^{7/8}\log t\cdot y^2\log y + \frac{t\log\log t\cdot
\log^2y}{\log t} \Big),
\end{align*}
the last step due to Lemma \ref{hp2.sum.error.lemma}. Since $r_1$ and
$r_2$ are distinct primes, the innermost sum is simply
$\pi(t;r_1r_2,1)$, and thus
\begin{multline}
\sum_{p\le t} h(p)^2~\le~\sum_{q_1,q_2\le y^2} \log q_1 \log q_2
\sum_{a_1,a_2\in\N} \sum_{\substack{r_1\in\P{q_1^{a_1}} \\
r_2\in\P{q_2^{a_2}}}} \pi(t;r_1r_2,1) \\
+ O\Big( t^{7/8}\log t\cdot y^2\log y + \frac{t\log\log t\cdot
\log^2y}{\log t} \Big).
\label{unhelpful.name}
\end{multline}

The contribution to the sum on the right-hand side of equation
\eqref{unhelpful.name} from those terms for which $q_1^{a_1}>t^{1/8}$
is
\begin{align*}
\sum_{q_1,q_2\le y^2} \log q_1 \log q_2 \sum_{\substack{a_1,a_2\in\N \\
q_1^{a_1}>t^{1/8}}} \sum_{\substack{r_1\in\P{q_1^{a_1}} \\
r_2\in\P{q_2^{a_2}}}} & \pi(t;r_1r_2,1) \\
&\ll~\sum_{q_1,q_2\le y^2} \log q_1 \log q_2
\sum_{\substack{a_1,a_2\in\N \\ q_1^{a_1}>t^{1/8}}}
\frac{t\log^2t}{q_1^{a_1}q_2^{a_2}} \\
&\ll~t\log^2t \sum_{q_1\le y^2} \sum_{\substack{a_1\in\N \\
q_1^{a_1}>t^{1/8}}} \frac{\log q_1}{q_1^{a_1}} \sum_{q_2\le y^2}
\sum_{a_2\in\N} \frac{\log q_2}{q_2^{a_2}} \\
&\ll~t\log^2t \sum_{q_1\le y^2} \frac{\log q_1}{t^{1/8}} \sum_{q_2\le
y^2} \frac{\log q_2}{q_2} \\
&\ll~t^{7/8}\log^2t \cdot y^2\log y
\end{align*}
by Lemma \ref{double.tricky.lemma}(a) and the estimates
\eqref{geometric.series}, \eqref{Chebyshev.bound}, and
\eqref{weaker.Chebyshev.bound}; the contribution from the terms for
which $q_1^{a_1}>t^{1/8}$ is bounded likewise. The remaining
contribution is
\begin{align*}
\sum_{q_1,q_2\le y^2} \log q_1 \log q_2 \sum_{\substack{a_1,a_2\in\N \\
q_1^{a_1},q_2^{a_2}\le t^{1/8}}} \sum_{\substack{r_1\in\P{q_1^{a_1}} \\
r_2\in\P{q_2^{a_2}}}} & \pi(t;r_1r_2,1) \\
&\ll~\sum_{q_1,q_2\le y^2} \log q_1 \log q_2
\sum_{\substack{a_1,a_2\in\N \\ q_1^{a_1},q_2^{a_2}\le t^{1/8}}}
\frac{t(\log\log t)^2}{q_1^{a_1}q_2^{a_2}\log t} \\
&\ll~\frac{t(\log\log t)^2}{\log t} \bigg( \sum_{q\le y^2}
\sum_{a\in\N} \frac{\log q}{q^a} \bigg)^2 \\
&\ll~\frac{t(\log\log t)^2\log^2 y}{\log t}
\end{align*}
by Lemma \ref{double.tricky.lemma}(b) and the estimates  
\eqref{geometric.series} and
\eqref{weaker.Chebyshev.bound}. Using both these bounds in equation
\eqref{unhelpful.name}, we conclude that
\[
\sum_{p\le t} h(p)^2~\ll~t^{7/8}\log t\cdot y^2\log y +
\frac{t(\log\log t)^2 \log^2y}{\log t}.
\]

We now evaluate $M_2(x)$ using partial summation. We have
\begin{align*}
M_2(x)~=~\sum_{p\le x} \frac{h(p)^2}p~&=~\sum_{p\le e^e} \frac{h(p)^2}p
+ \frac1x \sum_{e^e<p\le x} h(p)^2 + \int_{e^e}^x \frac{dt}{t^2}
\sum_{e^e<p\le t} h(p)^2 \\
&\ll~1 + \frac1x \cdot \frac{x(\log\log x)^2\log y}{\log x} \\
&\qquad{}+ \int_{e^e}^x \frac{dt}{t^2} \Big( t^{7/8}\log t\cdot y^2\log
y + \frac{t(\log\log t)^2 \log^2y}{\log t} \Big) \\
&\ll~\frac{y^2\log y}{\log x} + y^2\log y \int_{e^e}^x
\frac{\log t\, dt}{t^{9/8}} + \log^2y \int_{e^e}^x \frac{(\log\log  
t)^2}{t\log
t}\,dt.
\end{align*}
Evaluating these two integrals explicitly, we obtain
\[
M_2(x)~\ll~\frac{y^2\log y}{\log x} + y^2\log y + \log^2 y\cdot
(\log\log x)^3 \ll y^3\log^2y
\]
as claimed.
\end{proof}

\section{Normal number of cycles for the power generator}

If $(u,n)=1$, then the sequence $u^i\mod n$ for $i=1,2,\dots$ is purely
periodic.  We denote the length of the period by $\ord(u,n)$, which
of course is the multiplicative order of $u$ in $(\Z/n\Z)^\times$.
Even when $(u,n)>1$, the sequence $u^i\mod n$ is eventually periodic,
and we denote the length of the eventual cycle by $\oo(u,n)$.  So,  
letting $n_{(u)}$
denote the largest divisor of $n$ coprime to $u$, we have  
$\oo(u,n)=\ord(u,n_{(u)})$. For example, let $u=2,\,n=24$.  The  
sequence $u^i\mod n$ is $2,4,8,16,8,16,\dots$ with cycle length~2, and  
so $\oo(2,24)=\ord(2,3)=2$.

When iterating the $\ell$th power map modulo $n$, the length of the  
eventual cycle starting with $x=u$ is given by $\oo(\ell,\oo(u,n))$.   
We would like to have a criterion for when a residue is part of  
some cycle, that is, for when a residue is eventually sent back  
to itself when iterating $x\mapsto x^\ell\mod n$.

\begin{lemma}
A residue $u$ is part of some cycle under iteration of  
the map $x\mapsto x^\ell\mod n$ if and only if 
$(\ell,\oo(u,n))=1$ and, with $d=(u,n)$, we have
$(d,n/d)=1$.
\label{cycle.criterion.lemma}
\end{lemma}

\begin{proof}
If $(u,n)=d$, then high powers of $u$ will be $\equiv0\mod{n/n_{(d)}}$.
Thus, for $u$ to be in a cycle it is necessary that  $n/n_{(d)}=d$,
that is, $(d,n/d)=1$.  Further, it is necessary that  
$(\ell,\oo(u,n))=1$.
Indeed, if $\sigma=\oo(u,n)$, we would need $\ell^i\mod\sigma$ to be
purely periodic, which is equivalent to $(\ell,\sigma)=1$.  This proves
the necessity of the condition.  For the sufficiency, we have just
noted that $(\ell,\sigma)=1$ implies that $\ell^i\mod\sigma$ is
purely periodic.  This implies in turn that the sequence $u^{\ell^i}\mod{n_{(u)}}$
is purely periodic.  But the condition $(d,n/d)=1$ implies that $n_{(u)}=n/d$,
and as each $u^{\ell^i}\equiv0\mod d$, we have that $u^{\ell^i}\mod n$ is
purely periodic.
\end{proof}

For $d|n$ with $(d,n/d)=1$, let $C_d(\ell,n)$ denote the number of
cycles in the $\ell$th power map mod $n$ that involve residues $u$
with $(u,n)=d$. For the lower bound in Theorem  
\ref{number.of.cycles.theorem}
we shall deal only with $C_1(\ell,n)$, that is, cycles involving
numbers coprime to $n$.

\begin{lemma}
We have $C_1(\ell,n) \ge \phi(n)_{(\ell)}/\lambda(\lambda(n))$.
\label{c1.lower.bound.lemma}
\end{lemma}

\begin{proof}
It is easy to see that the subgroup of
$(\Z/n\Z)^\times$ of residues $u$ with $(\ell,\ord(u,n))=1$ has size
$\phi(n)_{(\ell)}$. (In fact, this is true for any finite abelian
group $G$: the size of the subgroup of elements with order coprime to
$\ell$ is $|G|_{(\ell)}$.) As the length of {\em any} cycle in
the $\ell$th power map is bounded above by $\lambda(\lambda(n))$, the  
lemma follows immediately.
\end{proof}

To investigate the normal size of $\phi(n)_{(\ell)}$, we introduce the  
function
\[
f_\ell(n)=\sum_{p\mid\ell} v_p(\phi(n))\log p.
\]
We also make use of the notation $q^a\| n$, which means that $q^a$ is the exact power of $q$ dividing $n$, that is, $q^a$ divides $n$ but $q^{a+1}$ does not.

\begin{prop}
For any fixed $\ell$, we have $f_\ell(n) \le (\log\log n)^2$ for  
almost all $n$, in fact for all but $O_\ell(x/\log\log x)$ integers $n\le x$.
\label{fn.small.aa.prop}
\end{prop}

\begin{proof}
We have
\begin{align*}
\sum_{n\le x}f_\ell(n)~&=~\sum_{p\mid\ell}\sum_{n\le x}\sum_{q^a\|  
n}v_p(\phi(q^a))\log p
~\le~x\sum_{p\mid\ell}\log p\sum_{q^a\le x}\frac{v_p(\phi(q^a))}{q^a}\\
&\le~x\sum_{p\mid\ell}\log p\sum_{p^a\le x}\frac{a-1}{p^a}+
x\sum_{p\mid\ell}\log p\sum_{q\le x}\frac{v_p(q-1)}{q}.
\end{align*}
Now
\[
x\sum_{p\mid\ell}\log p\sum_{p^a\le x}\frac{a-1}{p^a}
~\ll_\ell~x
\]
and, by \eqref{BT.prime.power},
\begin{align*}
x\sum_{p\mid\ell}\log p\sum_{q\le x}\frac{v_p(q-1)}{q}
&=~x\sum_{p\mid\ell}\log p\sum_{a\ge1}\sum_{q\in\P{p^a},\,q\le  
x}\frac1q\\
&\ll~x\sum_{p\mid\ell}\log p\sum_{a\ge1}\frac{\log\log x}{p^a}
~\ll_\ell~x\log\log x.
\end{align*}
Hence,
\[
\sum_{n\le x}f_\ell(n)~\ll_\ell~x\log\log x,
\]
so that the number of $n\le x$ with $f_\ell(n)>(\log\log n)^2$ is
$O_\ell(x/\log\log x)$.
\end{proof}

It is interesting that one can prove an Erd\H os--Kac theorem for  
$f_\ell(n)$
using as a tool the criterion of Kubilius--Shapiro (see \cite{K},  
\cite{S}).

\begin{proof}[Proof of the lower bound in Theorem  
\ref{number.of.cycles.theorem}]
Noting that $\phi(n)_{(\ell)}=\phi(n)/e^{f_\ell(n)}$, we have  
$\phi(n)_{(\ell)} \ge \phi(n)/\exp((\log\log n)^2)$ for almost all $n$  
by Proposition~\ref{fn.small.aa.prop}. Of course, $n\ge \phi(n) \gg  
n/\log\log n$ for all $n\ge3$.  
Therefore, using  Lemma~\ref{c1.lower.bound.lemma} and
Theorem~\ref{number.of.cycles.theorem}, we have
\begin{align*}
C(\ell,n)~\ge~C_1(\ell,n)~\ge~ 
\frac{\phi(n)_{(\ell)}}{\lambda(\lambda(n))}&\ge~
\frac{\phi(n)}{\exp((\log\log n)^2)\lambda(\lambda(n))}\\
&=~\frac{\phi(n)/n}{\exp((\log\log  n)^2)}\frac{n}{\lambda(\lambda(n))}\\
&=~\exp((1+o(1))(\log\log  n)^2\log\log\log n)
\end{align*}
for almost all $n$.  This
completes the proof of the lower bound in Theorem  
\ref{number.of.cycles.theorem}.
\end{proof}

We now consider the upper bounds in Theorem  
\ref{number.of.cycles.theorem}, first establishing a lemma.

\begin{lemma}
\label{fpsanalog}
Suppose $m$ is a positive integer and  $(d,m)=1$. For any integer $j\mid\lambda(m)$, the number of integers $u\in[1,m]$ with $(u,m)=1$ and $\ord(du,m)\mid\lambda(m)/j$ is at most $\phi(m)/j$.
\end{lemma}

\begin{proof}
In fact, we prove a more general statement for any finite abelian group $G$: let $\lambda(G)$ denote the exponent of $G$, that is, the order of the largest cyclic subgroup of $G$, or equivalently the least common multiple of the orders of the elements of $G$. Then for any $d\in G$ and any $j\mid\lambda(G)$, the number of elements $u\in G$ for which the order of $du$ divides $\lambda(G)/j$ is at most $\#G/j$. It is clear that the lemma follows immediately from this statement upon taking $G$ to be $(\Z/m\Z)^\times$. It is also clear that in this statement, the element $d$ plays no role whatsoever except to shuffle the elements of $G$ around, and so we assume without loss of generality that $d$ is the identity of~$G$.

Let $p$ be any prime dividing $\lambda(G)$, and choose $a\le b$ so that $p^a\|j$ and $p^b\|\lambda(G)$. When we write $G$ canonically as isomorphic to the direct product of cyclic groups of prime-power order, at least one of the factors must be isomorphic to $\Z/p^b\Z$. In every such factor, only one out of every $p^a$ elements has order dividing $\lambda(G)/j$, since all but $p^{b-a}$ elements of the factor have order divisible by $p^{b-a+1}$. Since there is at least one such factor for every $p^a\|j$, we conclude that at most one out of every $j$ elements of $G$ has order dividing $\lambda(G)/j$, as claimed.
\end{proof}

Note that this result in the case $d=1$ is Lemma 1 in~\cite{FPS}.  The above proof,
while similar in spirit to the proof in \cite{FPS}, is simpler.

Let $\tau(m)$ denote the number of positive divisors of $m$.

\begin{prop}
\label{cycleupperbd}
For any integers $\ell,n\ge2$ we have
$C(\ell,n)\le n\tau(\lambda(n))\tau(n)/\oo(\ell,\lambda(n))$.
\end{prop}

\begin{proof}
It is sufficient to show
that for each $\ell,n\ge2$ and each $d\mid n$ with $(d,n/d)=1$, we have
\begin{equation}
\label{cycled}
C_d(\ell,n)~\le~\frac{n\tau(\lambda(n))}{\oo(\ell,\lambda(n))}.
\end{equation}
Let $d\mid n$ with $(d,n/d)=1$.  We have seen in Lemma~\ref{cycle.criterion.lemma}
that for a residue $u\mod n$ with $(u,n)=d$ to be involved in a cycle,
it is necessary and sufficient that $(\ell,\ord(u,n/d))=1$.
For each integer $j\mid\lambda(n/d)$, let $C_{d,j}(\ell,n)$ denote the
number of cycles corresponding to residues $u$ with $(u,n)=d$ and
$\ord(u,n/d)=\lambda(n/d)/j$.  Writing such a residue $u$ as $du_1$,
we have $u_1\in[1,n/d]$ and $(u_1,n/d)=1$.  Thus, by Lemma  
\ref{fpsanalog},
we have that the number of such residues $u$ is at most $\phi(n/d)/j\le  
n/dj$.
Hence we have
\[
C_{d,j}(\ell,n)~\le~\frac{n/dj}{\ord(\ell,\lambda(n/d)/j)}.
\]
Now $\lambda(n/d)=\lambda(n)/d_1$ for some integer $d_1\le d$.  It is  
shown
in (15) of \cite{KP} that for $k\mid m$ we have  
$\oo(a,m/k)\ge\oo(a,m)/k$ for
any nonzero integer $a$.  Hence
\[
\ord(\ell,\lambda(n/d)/j)~=~\ord(\ell,\lambda(n)/ 
d_1j)~\ge~\oo(\ell,\lambda(n))/d_1j,
\]
so that
\[
C_{d,j}(\ell,n)~\le~\frac{n/dj}{\oo(\ell,\lambda(n))/ 
d_1j}~\le~\frac{n}{\oo(\ell,\lambda(n))}.
\]
Letting $j$ range over all divisors of $\lambda(n/d)$, we get that
\[
C_d(\ell,n)~\le~\frac{n\tau(\lambda(n/d))}{\oo(\ell,\lambda(n))},
\]
which immediately gives \eqref{cycled}.
\end{proof}

\begin{proof}[Proof of the upper bounds in Theorem  
\ref{number.of.cycles.theorem}]
Note that from
\cite[Theorem 4.1]{EP}, we have $\tau(\lambda(n))<\exp((\log\log n)^2)$
for almost all $n$.  Furthermore, letting $\Omega(n)$ denote the number  
of prime factors of $n$ counted with
multiplicity, we know that the normal order of $\Omega(n)$ is $\log\log  
n$; in particular, we have $\Omega(n) < \log\log n/\log 2$ for almost  
all $n$. Since the inequality $\tau(n) \le 2^{\Omega(n)}$ is  
elementary, this implies that $\tau(n)<\log n$ for almost all $n$. We  
conclude from Proposition~\ref{cycleupperbd} that
\[
C(\ell,n)~<~n\exp(2(\log\log n)^2) / \oo(\ell,\lambda(n))
\]
for almost all~$n$.

The three upper bounds in
Theorem~\ref{number.of.cycles.theorem} therefore follow respectively from
three results in the new paper of Kurlberg and the second  
author \cite{KP}:
Theorem 4 (1), which states that for any function $\ep(n)\to0$, we have
$\oo(\ell,\lambda(n))\ge n^{1/2+\ep(n)}$ almost always;
Theorem 22, which states that a positive proportion of integers $n$
have $\oo(\ell,\lambda(n))\ge n^{.592}$; and
Theorem 28, which states that if the GRH is true, then
\[
\oo(\ell,\lambda(n))~=~n/\exp((1+o(1))(\log\log n)^2\log\log\log n)
\]
on a set of asymptotic density~1.  (Note that the proof of this result uses Theorem
\ref{lambda.lambda.normal.order.theorem} of the current paper.)
\end{proof}

\section{Higher iterates}

Here we sketch what we believe to be a viable strategy for establishing  
an analogue of Theorem~\ref{lambda.lambda.normal.order.theorem} for the  
higher iterates $\lambda_k$ where $k\ge3$.
As in the case of $k=2$, we have generally that
\[
\frac{n}{\lambda_k(n)}~=~\frac{n}{\phi_k(n)}\frac{\phi_k(n)}{\lambda_k(n 
)}.
\]
We always have $n/\phi_k(n)\le(c\log\log n)^k$, which is already a good  
enough estimate for our purposes. Even better, however, it is known  
\cite{EGPS} that for each fixed $k$, we have  
$n/\phi_k(n)\ll(\log\log\log n)^k$
for almost all $n$.  The problem therefore reduces to comparing  
$\lambda_k(n)$ to $\phi_k(n)$.
Probably it is not hard to get analogs of Propositions  
\ref{same.large.prime.divisors.prop}
and \ref{large.primes.in.phi.phi.prop}, where we replace $y^2$ with  
$y^k$.  The problem
comes in with the proliferation of cases needed to deal with small  
prime factors.
As with the second iterate, we expect the main contribution to come  
from the
``supersquarefree" case.  In particular, let
\[
h_k(n) ~=~ \sum_{p_1\mid n} \sum_{p_2\mid p_1-1} \dots \sum_{p_k\mid
p_{k-1}-1} \sum_{q\le y^k} v_q(p_k-1) \log q.
\]
We expect $h_k(n)$ to be the dominant contribution to  
$\log(\phi_k(n)/\lambda_k(n))$
almost always.  But it seems hard not only to prove this in general but  
also to
establish the normal order of $h_k(n)$.

It would seem useful in this endeavor to have a uniform estimate of the  
shape
\begin{equation}
\label{strongrecipsum}
\sum_{p\in\P{m},\,p\le x}\frac1p~\sim~\frac{\log\log x-\log\log  
m}{\phi(m)}
~\text{ for }~x\ge m^{1+\ep}.
\end{equation}
Even under the assumption of the Riemann Hypothesis for Dirichlet $L$-functions, 
\eqref{strongrecipsum} seems difficult, and maybe it  
is false.
It implies with $x=m^2$ that the sum is $\ll1/\phi(m)$,  
when all
we seem to be able to prove, via sieve methods, is that it is  
$\ll(\log\log m)/\phi(m)$.

Assuming uniformity in \eqref{strongrecipsum}, it seems that on average
\[
h_k(n) \sim \frac1{(k-1)!}(\log\log n)^k\log\log\log n,
\]
supporting Conjecture \ref{lambdak.normal.order.conj}. It would be a  
worthwhile
enterprise to try to verify or disprove the Conjecture in the case
$k=3$, which may be tractable.

Going out even further on a limb, it may be instructive to think of what
Conjecture \ref{lambdak.normal.order.conj} has to say about the normal  
order
of $L(n)$, the minimum value of $k$ with $\lambda_k(n)=1$.  The  
expression
$(1/(k-1)!)(\log\log n)^k\log\log\log n$ reaches its maximum value when
$k\approx\log\log n$.  Is this formula then trying to tell us that
we have $L(n)\ll \log\log n$ almost always?  Perhaps so.

There is a second argument supporting the thought that $L(n)\ll\log\log  
n$
almost always.  Let $P(n)$ denote the largest prime factor of an  
integer $n>1$, and let $\ell(n)=P(n)-1$ for $n>1$, $\ell(1)=1$. 
Clearly,
$\ell(n)\mid\lambda(n)$ for all $n$, so that if $L_0(n)$ is the least
$k$ with $\ell_k(n)=1$, then $L_0(n)\le L(n)$.  It may be that the
difference $L(n)-L_0(n)$ is usually not large.  In any event, it seems
safe to conjecture that $L_0(n)$ is usually of order of magnitude  
$\log\log n$, due to the following argument.
For an odd prime $p$, consider the quantity $\log\ell(p)/\log
p\approx\log P(p-1)/\log (p-1)$. It may be that this quantity is
distributed as $p$ varies through the primes
in the same way that $\log P(n)/\log n$ is distributed as $n$ varies
through the
integers, namely the Dickman distribution.
Such a conjecture has been made in various papers.
If so, it may be that the sequence
\[
\frac{\log\ell(p)}{\log p},~
\frac{\log\ell_2(p)}{\log\ell(p)},~\dots
\]
behaves like a sequence of
independent random variables, each with the Dickman distribution.
And if so, it may then be reasonable to assume that almost always
we get down to small numbers and terminate in about $\log\log n$ steps.
A similar probabilistic model is considered in \cite{B}, but for the
simpler experiment of finding the joint distribution of logarithmic
sizes of the various prime factors of a given number $n$.

At the very least, we can prove that $L(n)\ll\log\log n$ infinitely  
often.

\begin{proof}[Proof of Theorem \ref{log.log.iterates.theorem}]
Notice that the definition of $\lambda(n)$ as a least common multiple,
together with the fact that $\lambda(p^a)\mid \lambda(p^{a+1})$ always,
implies that
\[
\lambda\big( \lcm\{ m_1,\dots,m_j \} \big)~=~\lcm\big\{ \lambda(m_1),
\dots, \lambda(m_j) \big\}
\]
for any positive integers $m_1,\dots,m_j$. A trivial induction then
shows that
\[
\lambda_k\big( \lcm\{ m_1,\dots,m_j \} \big)~=~\lcm\big\{
\lambda_k(m_1), \dots, \lambda_k(m_j) \big\}
\]
for any $k\ge0$. Since the least common multiple of a set of numbers
equals 1 precisely when each number in the set equals 1, we deduce that
\[
L\big( \lcm\{ m_1,\dots,m_j \} \big)~=~\max\big\{ L(m_1), \dots, L(m_j)
\big\}.
\]

We apply this identity with $m_i=i$.  
Let  $n_j = \lcm\{1,2,\dots,j\}$. We have $\log n_j = \sum_{i\le j}
\Lambda(i)$, which is asymptotic to $j$ by the prime number theorem.
On the other hand, it is trivial
that for any number $n$ we have $L(n) \le1+ (1/\log 2)\log n$, as
$\lambda_{i+1}(n)\le(1/2)\lambda_i(n)$ for $1\le i<L(n)$.  Therefore
\begin{align*}
L(n_j)~&=~\max\{ L(1),\dots,L(j) \}~\le~
1+\max\left\{\frac{\log1}{\log2},\dots,\frac{\log j}{\log2}\right\}\\
&=~1+\frac{\log j}{\log2}
~=~\left(\frac1{\log2}+o(1)\right)\log\log n_j.
\end{align*}
\end{proof}

We can improve on the estimate in Theorem  
\ref{log.log.iterates.theorem},
but not by much.  Say we let $N_j$ be the product of all primes $p\le  
j^{3.29}$
with $p-1\mid n_j$, with $n_j$ as in the above proof.  It  
follows
from Friedlander \cite{F} that a positive proportion of the primes  
$p\le j^{3.29}$
have the required property.  Thus, $N_j>\exp(cj^{3.29})$ for some positive
constant $c$ and all sufficiently large values of $j$.  But  
$\lambda(N_j)\mid n_j$,
so that $L(N_j)\le 2+j/\log 2$.  Hence  
$L(N_j)<.439\log\log N_j$
for $j$ sufficiently large.  (This result can be improved by a very  
small
margin using a more recent result of Baker and Harman \cite{BH}, but the
argument is a bit more difficult, since they do not get a positive
proportion of the primes with the required property.) It is likely that  
$L(n)\ll\log\log\log n$ infinitely often, possibly even that $L(n)  
\ll_k \log_k n$ infinitely often for arbitrary $k$-fold-iterated  
logarithms.

One may also study the maximal order of $L(n)$.  The analogous problem  
for
the iterated $\phi$-function is relatively trivial, but not so for  
$\lambda$.
If there can exist very long ``Sophie Germain chains", that is,  
sequences
of primes $p_1,p_2,\dots,p_k$ where each $p_i=2p_{i-1}+1$, for $i>1$,  
then
we might have $L(p_k)\sim(1/\log 2)\log p_k$.  We might even perturb  
such
a chain by a small amount and keep the asymptotic relation, say by  
occasionally
having $p_i=4p_{i-1}+1$.  It seems hard to prove that long enough chains
to get the the asymptotic for $L(p_k)$ do not exist, but probably they
don't on probabilistic grounds.  We can at least say that
$L(n)\ge1+(1/\log 3)\log n$ infinitely often, since this inequality is  
attained
when $n$ is a power of 3.

\end{document}